\documentclass[journal]{IEEEtran}

\usepackage{soul}
\usepackage{amssymb}
\usepackage{hyperref}
\usepackage{gensymb}
\usepackage{graphicx}
\usepackage{caption}
\usepackage{subfig}
\usepackage{tabularx}
\usepackage{color}
\usepackage{xcolor}
\usepackage{amsmath}
\usepackage[font={small}]{caption}

\ifCLASSINFOpdf
\else
\fi
\begin{document}

\title{\LARGE Uncertainty Sets For Wind Power Generation}
\author{\normalsize Yury Dvorkin, \textit{Student Member, IEEE}, Miles Lubin, Scott Backhaus, Michael Chertkov, \textit{Senior Member, IEEE}. \vspace{-30pt}}

\markboth{ Submitted to IEEE Power and Engineering Letters} 
{Shell \MakeLowercase{\textit{et al.}}: Bare Demo of IEEEtran.cls for Journals} 

\maketitle 

\begin{abstract}
As penetration of wind power generation increases, system operators must account for its stochastic nature in a reliable and cost-efficient manner. These conflicting objectives can be traded-off by accounting for the variability and uncertainty of wind power generation. This letter presents a new methodology to estimate uncertainty sets for parameters of probability distributions that capture wind generation uncertainty and variability.
\vspace{-10pt}
\end{abstract}

\begin{IEEEkeywords}
Wind power uncertainty, wind power variability, power system operations. \vspace{-10pt}
\end{IEEEkeywords}

\IEEEpeerreviewmaketitle
\section{Introduction}

\IEEEPARstart{W}{ind} power generation (WPG) introduces variability and uncertainty in operational planning. As defined in \cite{ela_2012}, \textit{variability} of WPG is the random fluctuation of wind speed caused by physical processes in the atmosphere while \textit{uncertainty} of WPG results from wind forecast errors. To account for variability and uncertainty, approaches to robust unit commitment (RUC) \cite{guan_2014}, chance constrained optimal power flow (CC-OPF) \cite{bienstock_2009} and distributionally robust CC-OPF (RCC-OPF) \cite{bienstock_2009},\cite{new_submission} have been formulated and tested. However, the performance of these models depend on the accuracy of parameters of probability distributions that define the uncertainty sets used; therefore, the uncertainty in distribution parameters must be accounted for. \textcolor{black}{This letter makes two contributions: i) we use a data-driven analysis to relate intra-hour wind speed variability to the hourly-average wind speed and ii) using this relationship we construct uncertainty sets that can be interpreted in terms of variability and uncertainty of WPG.} 

Bienstock et al \cite{bienstock_2009} show that the CC-OPF formulation with Gaussian-distributed deviations of WPG and precisely known mean and variance can be extended to the RCC-OPF formulation, where the parameters of the Gaussian deviations (both mean and variance) fall within uncertainty sets. This RCC-OPF is implemented and tested on a large-scale system in \cite{new_submission}. \textcolor{black}{In this letter, we interpret the set for the mean in terms of the wind power uncertainty, while the set for the variance is deemed as the wind power variability. This approach differs from the uncertainty sets on a random variable in \cite{guan_2014} by introducing uncertainty sets on distribution parameters describing the random variable. We derive these sets with a data-driven approach such that the resulting RCC-OPF remains tractable \cite{bienstock_2009, new_submission}, yielding a practical approach for modeling WPG. Although the RCC-OPF model is nominally more conservative than the CC-OPF model, its solution incurs a lower operating cost when tested against actual realizations of WPG \cite{new_submission}.}

\vspace{-5pt}
\section{Methodology}

\subsubsection{Data}
\label{sec::data}
{We use historical wind speed measurements at 5-minute resolution from the Goodnoe meteorological station in the Bonneville Power Authority (BPA) system \cite{bpa} and one-hour resolution wind speed forecasts produced by the NOAA Rapid Refresh numerical weather prediction model \cite{noaa} for the same location. The historical measurements and forecasts are detrended by using data from the same calendar season (December--February) and range of hours (00:00--04:00 AM). }

\subsubsection{Wind Speed Variability}
\label{sec::variability}
\begin{figure}[b]
\vspace{-20pt}
\begin{tabular}{c}
\subfloat{\includegraphics[width=0.95\linewidth]{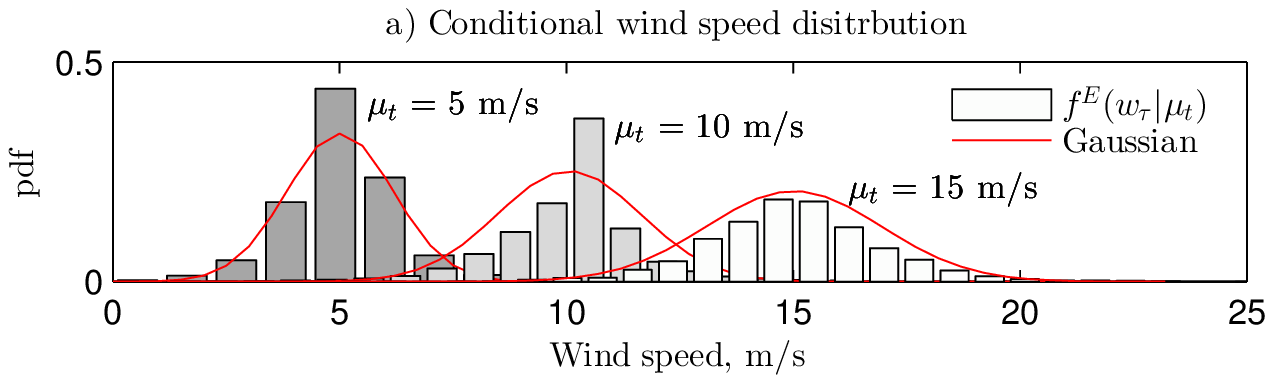}}\\
\subfloat{\includegraphics[width=0.95\linewidth]{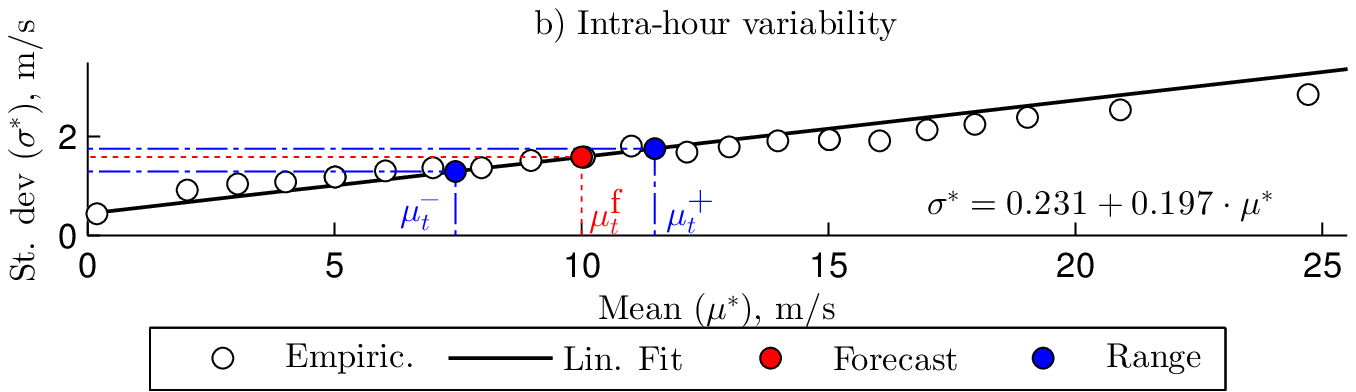}}\\
\end{tabular}
\caption{a) Empirical distribution $f^E(w_{t} | \mu^*)$ (histograms) and their Gaussian best fits for a few $\mu^*$. b) Empirical relationship between $\sigma*$ and $\mu^*$ and its linear fit $\sigma^*(\mu^*)$.   }
\label{fig:BPApdfs}
\end{figure}

{For any subinterval $\tau$ of hour-long interval $t$ wind speed $w_{\tau}$ can be written as $w_{\tau}= \mu_t + \epsilon_{\tau}$ \cite{bienstock_2009}, where $\mu_t$ and $\epsilon_{\tau}$ are the hourly-average wind speed and intra-hour, zero-mean variability around $\mu_t$ \cite{ela_2012}, respectively. \textcolor{black}{To parametrize $w_{\tau}$, we calculate the average wind speed ($\mu_t$) for each hourly interval within the studied period using the 5-minute resolution Goodnoe data \cite{bpa}. Next, we bin the hourly intervals based on their hourly-averaged wind speed $\mu_t$. This binning is in the range from 0 m/s to 25 m/s (the  typical operating range of wind turbines) and has the width of 1 m/s. For each bin we obtain a histogram that represents an empirical intra-hour wind speed distribution $f^E(w_{\tau}|\mu_t)$ conditioned by the hourly-average wind speed ($\mu_t$).  Each histogram is then fit to the Gaussian distribution, yielding  estimated mean ($\mu^*$) and standard deviation values ($\sigma^*$) for every bin. Both the histogram and its Gaussian distribution are shown in Fig.~\ref{fig:BPApdfs}a) for several values of $\mu_t$. Figure~\ref{fig:BPApdfs}b) shows that $\sigma^*$ scales linearly with $\mu^*$ as $\sigma^*(\mu^*) = 0.231 + 0.197 \cdot \mu^*$ for $\mu^* \in \left[0, 25\right]$  m/s, which is consistent with velocity distributions for high Reynolds number atmospheric flows \cite{kaimal}.} If there is no error in a given hourly-averaged wind speed forecast $\mu_t^f$, the wind speed variability is estimated using the linear mapping $\sigma^*(\mu^*)$ shown in Fig.~\ref{fig:BPApdfs}b). The resulting distribution is then given as $N [w_{\tau}; \mu_t^f, (\sigma^*(\mu_t^f))^2 ]$, and can be used in \cite{guan_2014} and \cite{bienstock_2009} for uncertainty sets accounting for wind speed variability.

\subsubsection{Wind Speed Uncertainty and Variability}
\label{sec::uncertainty}
{Wind forecast errors cause $\mu_t\neq \mu_t^f$. The wind forecast error $e_t$ is calculated from the NOAA \cite{noaa} data as $e_{t}(\Delta T)=\mu_t-\mu_t^f(\Delta T)$, where $\mu_t^f(\Delta T)$ is the forecast for hour $t$ made $\Delta T$ hours in advance. The empirical distribution $f^E \left(e_{t};\Delta T \right)$ for $\Delta T=1$ hour is shown in Fig.~\ref{fig:NOAApdds}a). We propose to represent $f^E(\cdot)$ with a generalized Gaussian distribution, $f^G(\cdot)$:}
\begin{equation}
\label{eq1}
f^G(e_t;\Delta T,\mu^-,\mu^+)=\frac{  \int_{\mu^-}^{\mu^+}d\mu   N[e_t;\mu,\sigma^*(\mu)]  }{(\mu^+ - \mu^-)},
\end{equation}
{\noindent where $\sigma^*(\mu)$ is the fit from Fig.~\ref{fig:BPApdfs}b) and the dependence of $\mu^+$ and $\mu^-$ on $\Delta T$ is suppressed. \textcolor{black}{The advantage of $f^G(\cdot)$ over a single Gaussian distribution is that it improves the goodness-of-fit to empirical data \cite{hahn_1967} by encapsulating a linear superposition of Gaussian distributions over a range of means $[\mu^-, \mu^+]$ that represents the wind uncertainty. Furthermore,  \cite{hahn_1967} explains that other single non-Gaussian distributions can be modelled using $f^G(\cdot)$.} The best fit $\mu^-$ and $\mu^+$ are computed by solving the optimization problem:}
\begin{equation}
\label{eq2}
\arg\min_{\!\!\!\!\!\!\mu^-, \mu^+} \int  [f^G(e_t;\Delta T,\mu^-, \mu^+) -f^E(e_t;\Delta T)]^2 de_t,
\end{equation}
{which minimizes the mean square difference between $f^E(\cdot)$ and $f^G(\cdot)$ and ensures a better fit to the historical data than a single Gaussian distribution, as shown in  Fig. \ref{fig:NOAApdds}a). We interpret the range $[\mu_t^{f-},\mu_t^{f+}]$, where $\mu_t^{f-}=\mu^f_t+\mu^-$ and $\mu_t^{f+}=\mu^f_t+\mu^+$, as the bounds of the uncertainty set for the mean wind speed. Fig.~\ref{fig:NOAApdds}b) displays $\mu_t^{f-}$ and $\mu_t^{f+}$ for $\mu^f_{t} = 10$ m/s for different $\Delta T$. Using $\sigma^*(\mu^*)$ from Fig.~\ref{fig:BPApdfs}b), we compute the bounds on $\sigma^*$ as $\sigma^{*} (\mu_t^{f-})$ and $\sigma^{*} (\mu_t^{f+})$, which are shown in Fig.~\ref{fig:BPApdfs}c). The  ranges $[\mu_t^{f-},\mu_t^{f+}]$ and $[\sigma^{*} (\mu_t^{f-}),\sigma^{*} (\mu_t^{f+})]$, if converted to wind generation as explained below, can be used in the RCC-OPF formulation from \cite{bienstock_2009}, \cite{new_submission}.
\begin{figure}[t]
\begin{tabular}{c}
\subfloat{\includegraphics[width=0.95\linewidth]{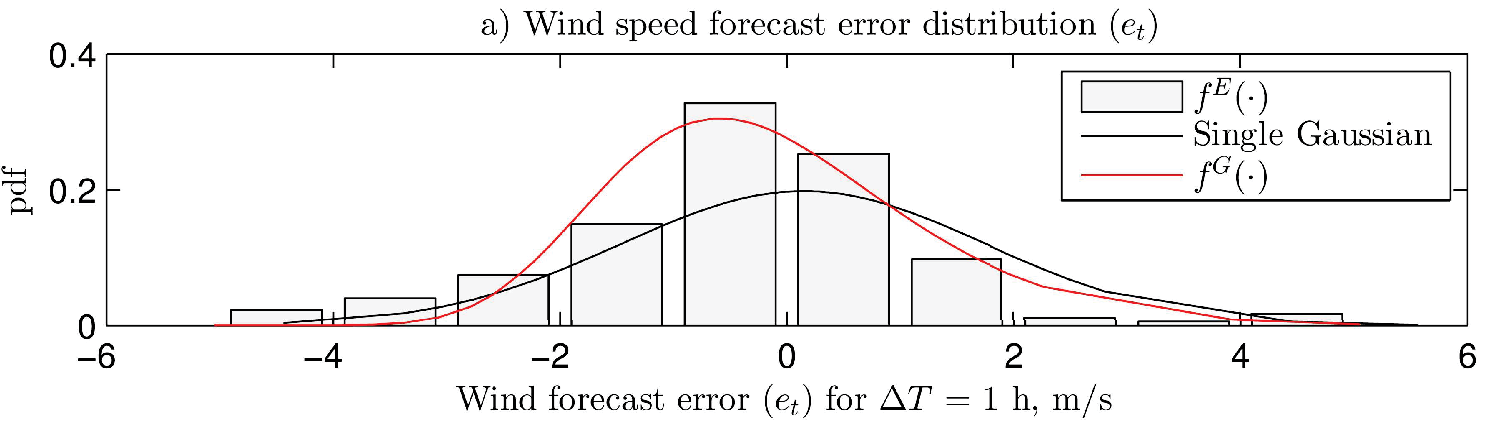}}\\
\subfloat{\includegraphics[width=0.95\linewidth]{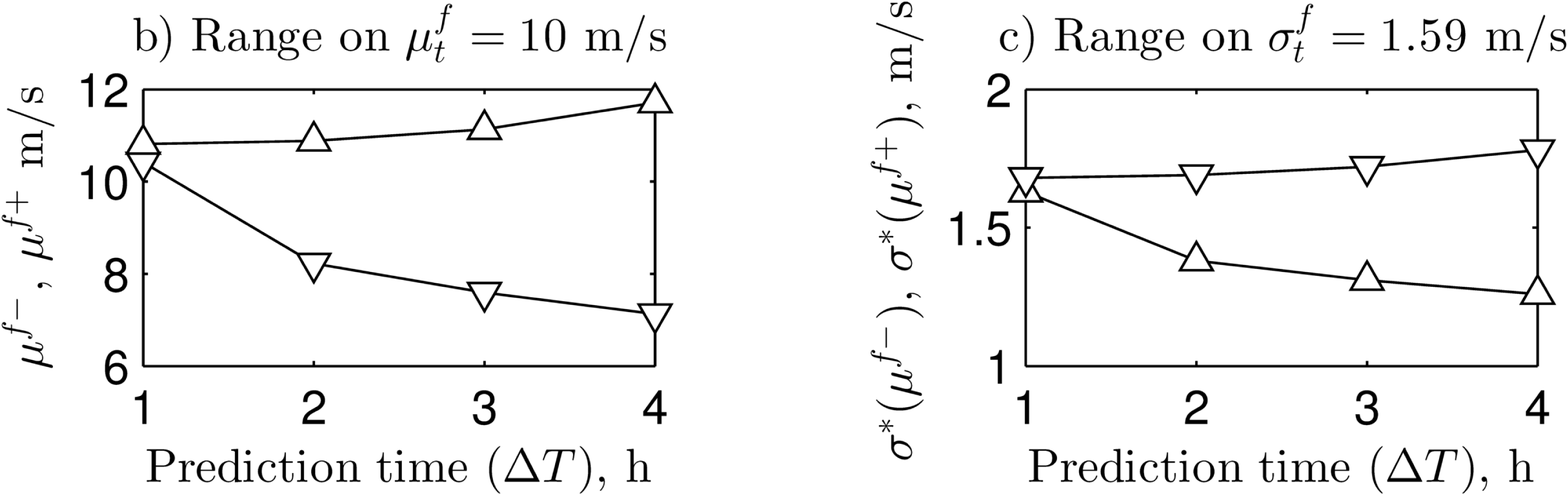}}\\
\end{tabular}
\caption{a) Empirical distribution $f^E(e_t;\Delta T)$ (histogram), the best fit single Gaussian distribution (black) and the best fit distribution $f^G(\cdot)$ (red) from Eq.~(\ref{eq1}) for $\Delta T$= 1 h. b) $\mu_t^{f-}$ and $\mu_t^{f+}$ versus $\Delta T$ for $\mu^f_{t} = 10$ m/s. c) $\sigma^*(\mu_t^{f-})$ and $\sigma^*(\mu_t^{f+})$ for the data in b).}
\label{fig:NOAApdds}
\end{figure}
\begin{figure}[t]
\vspace{-10pt}
\begin{tabular}{c}
\subfloat{\includegraphics[width=\linewidth]{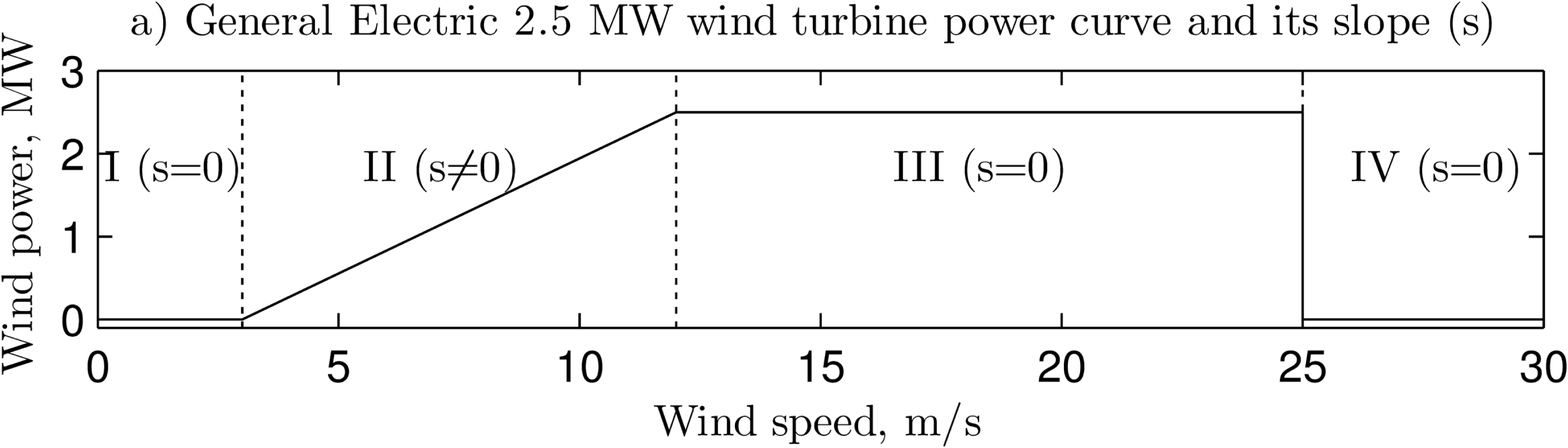}} \\
\subfloat{\includegraphics[width=\linewidth]{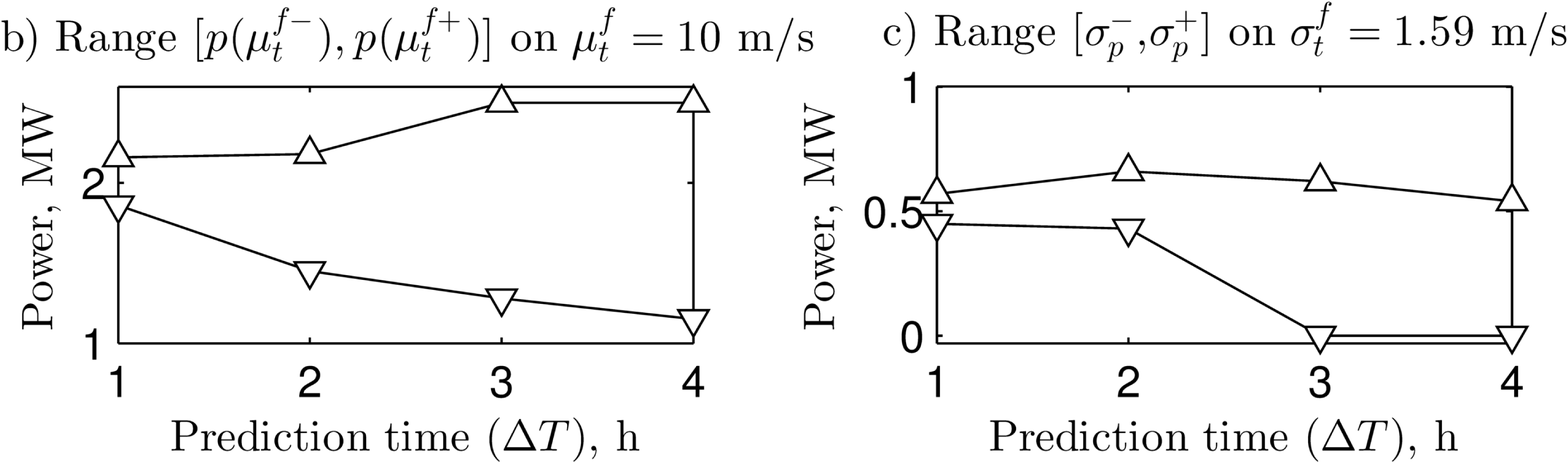}}\\
\end{tabular}
\caption{a) Typical wind power turbine curve $p(\mu)$ with four nominal operating regions \cite{ge}. b) Ranges for the on the mean WPG (i.e. forecast error) computed using Eq. (\ref{eq:power_interval}) for the data in Fig.~\ref{fig:NOAApdds}c). c) Ranges on WPG standard deviation (variability) computed using Eqs.~(\ref{eq:uset_sigma-}) and (\ref{eq:uset_sigma+}) for the data in b).    }
\label{fig:ranges}
\vspace{-5pt}
\end{figure}

\subsubsection{Conversion to Wind Power}
\label{sec::conversion}
We illustrate the conversion to wind power using the single wind turbine power curve $p(\mu)$ shown in Fig.~\ref{fig:ranges}a), \cite{ge}. This procedure can be generalized to multiple turbines by using an aggregated wind power curve \cite{hayes_2011}. The conversion of the range $[\mu_t^{f-},\mu_t^{f+}]$ is given by:
\begin{align}
[\mu_t^{f-},\mu_t^{f+}]\rightarrow [p(\mu_t^{f-}),p(\mu_t^{f+})].\label{eq:power_interval}
\end{align}
Figure~\ref{fig:ranges}b) shows $p(\mu_t^{f-})$ and $p(\mu_t^{f+})$ corresponding to $\mu_t^{f-}$ and $\mu_t^{f+}$, respectively, from Fig.~\ref{fig:NOAApdds}b). In Fig.~\ref{fig:NOAApdds}b), the growth of $p(\mu_t^{f+})$ at larger $\Delta T$ is eventually clipped by $p(w)$ as $\mu_t^{f+}$ enters Region III of the turbine curve.  If the entire range $[\mu_t^{f-},\mu_t^{f+}]$ is in Region III, then the range $[p(\mu_t^{f-}),p(\mu_t^{f+})]$ collapses to zero width around the maximum output.

The conversion of the range $[\sigma^{*} (\mu_t^{f-}),\sigma^{*} (\mu_t^{f+})]$ is shaped by the RCC-OPF formulation in \cite{new_submission}. Note that $\sigma_p$ can be obtained  $\forall p (\mu_t) \in [p(\mu_t^{f-}),p(\mu_t^{f+})]$ using the relation $\sigma^*(\mu)$ from Fig.~\ref{fig:BPApdfs}b) and the slope $s$ of the turbine curve as $\sigma_p=s (\mu_t) \cdot \sigma^*(\mu_t)$.  However, using the wind turbine power curve from Fig.~\ref{fig:ranges}a) results in difficult nonconvexity in distributionally robust formulations (e.g., Eq. (27)-(28) in  \cite{new_submission}). To avoid this nonconvexity, we assume that the range on the mean value in (Eq.~\ref{eq:power_interval}) and  range on the standard deviation are independent. Therefore, the range on the standard deviations is given by:
\begin{gather}
\sigma_p^- = \min_{\mu \in [\mu_t^{f-},\mu_t^{f+}]} s(\mu) \sigma^*(\mu )\label{eq:uset_sigma-}  \\
\sigma_p^+ = \max_{\mu \in [\mu_t^{f-},\mu_t^{f+}]} s(\mu ) \sigma^*(\mu ) \label{eq:uset_sigma+}
\end{gather}
Figure~\ref{fig:ranges}c) shows the conversion to the robust interval on the standard deviation of WPG for the data in Fig.~\ref{fig:NOAApdds}c).
\vspace{-10pt}
\section{Conclusion} We have presented a data-driven method to develop robust intervals for distribution parameters, which preserves the physical relationship between instantaneous and hour-average wind speed and power. This method is suitable for uncertainty sets in the RUC \cite{guan_2014} and RCC-OPF \cite{bienstock_2009}, \cite{new_submission}, which can be interpreted in terms variability and uncertainty of WPG. \textcolor{black}{The case study in \cite{new_submission} shows that the RCC-OPF model with these uncertainty sets outperforms several benchmarks in terms of several cost and reliability metrics.} 
\vspace{-10pt}
\section*{Acknowledgment}

This work was supported by the Advanced Grid Modeling Program in the U.S. Department of Energy Office of Electricity under Contract No. DE-AC52-06NA25396.
\vspace{-10pt}
\ifCLASSOPTIONcaptionsoff
  \newpage
\fi

\end{document}